\documentclass{article}

\usepackage{amssymb}
\usepackage{longtable}

\newcommand\N{{\mathbb N}}

\def\maxgenus{50}
\renewcommand{\leq}{\leqslant}
\renewcommand{\geq}{\geqslant}

\author{Maria Bras-Amor\'os
\thanks{Universitat Rovira i Virgili, Catalonia, Spain,
e-mail: maria.bras@urv.com}}

\date{Published at \\Semigroup Forum, Springer, vol. 76, n. 2, pp. 379-384, March 2008}

\title{Fibonacci-like behavior of the number of numerical semigroups of a given genus}

\begin{document}
\maketitle

\begin{abstract}
We conjecture a Fibonacci-like property on the number of numerical
semigroups of a given genus. Moreover we conjecture that the associated quotient
sequence approaches the golden ratio.
The conjecture is motivated by the
results on the number of semigroups of genus at most $\maxgenus$.
The Wilf conjecture has also been checked for all numerical
semigroups with genus in the same range.
\end{abstract}

\section{Introduction}

Let $\N_0$ denote the set of all non-negative integers.
A {\it numerical semigroup} is a subset $\Lambda$ of $\N_0$ containing $0$,
closed under addition and with finite complement in $\N_0$.
For a numerical semigroup $\Lambda$ define the {\it genus} of $\Lambda$ as the number
$g=\#(\N_0\setminus\Lambda)$.
As an example of numerical semigroup one can take
$$\{0,4,5,8,9,10\}\cup\{i\in\N_0:i\geq 12\}.$$
In this case the genus is $6$.

We are interested on the number $n_g$ of numerical semigoups of genus
$g$.
It is obvious that $n_0=1$ since ${\mathbb N}_0$ is the unique
numerical semigroup of genus $0$.
On the other hand, if $1$ is in a numerical semigroup, then any
non-negative integer must belong also
to the numerical semigroup, because any non-negative
integer is a finite sum of $1$'s. Thus, the unique numerical
semigroup with genus $1$ is $\{0\}\cup\{i\in\N_0:i\geq 2\}$
and $n_1=1$.
We conjecture
\begin{enumerate}
\item $n_g\geq n_{g-1}+n_{g-2}$, for $g\geq 2$
\item $\mbox{lim}_{g\rightarrow\infty}\frac{n_{g-1}+n_{g-2}}{n_g}=1$
\item $\mbox{lim}_{g\rightarrow\infty}\frac{n_{g}}{n_{g-1}}=\phi$, where
  $\phi$ is the golden ratio.
\end{enumerate}
Notice that point 2 would mean the sequence $n_g$
to behave asymptotically as the Fibonacci sequence.
This conjecture is motivated by the
results on the number of semigroups of genus at most $\maxgenus$.

\section{Computation of $n_g$}

Every numerical semigroup can be generated by a finite set of
elements and a minimal set of generators is unique (see for instance \cite{FrGoHa}).
Let $m$ be the smallest non-zero element of a numerical semigroup
$\Lambda$,
usually referred as its {\it multiplicity}.
The set of minimal generators of $\Lambda$
can be easily computed because
it is a subset of the finite set
$$\{m\}\cup\{\lambda\in\Lambda\setminus\{0\}:\lambda-m\not\in\Lambda\}.$$
This last set is intimately related to the
Ap\'ery set of $\Lambda$
\cite{Apery,FrGoHa,RoGaGaBr2002,RoGaGaBr2005,MaHe}.

The {\it conductor} of a numerical semigroup $\Lambda$
is the unique integer $c$
with $c-1\not\in\Lambda$ and $c+\N_0\subseteq\Lambda$.
Given a numerical semigroup $\Lambda$ of genus $g$ and conductor $c$,
$\Lambda\cup\{c-1\}$ is a numerical semigroup and its genus is $g-1$.
So, any numerical semigroup of genus $g$ can be obtained from
a numerical semigroup of genus $g-1$ by removing one element
larger than or equal to its conductor. It is easy to check
that when removing such an element from a numerical semigroup,
the set obtained is a numerical semigroup if and only if the removed element
belongs to the set of minimal generators.

For instance, the unique numerical semigroup of genus $0$ is $\N_0$.
Its unique minimal generator is $1$. Now, removing $1$ from $\N_0$ we obtain
$$\Lambda_1=\{0\}\cup\{i\in\N_0:i\geq 2\}.$$
The numerical semigroup $\Lambda_1$ is the unique numerical semigroup
of genus $1$.
In turn, $\Lambda_1$ has minimal generators $2$ and $3$.
By removing $2$ from $\Lambda_1$ we obtain the numerical semigroup
$$\Lambda_{2,1}=\{0\}\cup\{i\in\N_0:i\geq 3\}$$
and
by removing $3$ from $\Lambda_1$ we obtain the numerical semigroup
$$\Lambda_{2,2}=\{0,2\}\cup\{i\in\N_0:i\geq 4\}.$$
The semigroups $\Lambda_{2,1}$ and $\Lambda_{2,2}$
are all the numerical semigroups of genus $2$.

For the results in this paper we computed $n_g$ by brute approach.
That is, we generated all
numerical semigroups of genus $g$ from all numerical semigroups of
genus $g-1$ as explained and then we counted them.
This way one should be able to generate all semigroups of any given genus.

In the web page
\mbox{\small{\tt http://w3.impa.br/\~{}nivaldo/algebra/semigroups/index.html}}
by Nivaldo Medeiros one can find all numerical semigroups of genus
up to $12$.
In Neil Sloane's On-line Encyclopedia of Integer
Sequences \cite{Sloane:OnlineEncyclopediaIntegerSequences}
there are the values of $n_g$ for $g\leq 14$.
We could compute all numerical semigroups of genus up to $\maxgenus$.

The obstruction on the calculus of all numerical semigroups
of a given {\it large} genus using the method explained above
is the huge size of the results and the
need to keep them for the next step. Indeed,
the growth of $n_g$ is apparently exponential with $g$ and doing
computations beyond a certain genus
is really difficult using the
computational means one can find at present.
For instance, using a Pentium D 3.00 GHz
with 1 GB of RAM it took $19$ days to compute
all semigroups of genus $50$ and we expect that it would take about one month
to compute all numerical semigroups of genus $51$.
The size of the compressed file containing all numerical semigroups of
genus $50$ is $3.6$ GB.

\begin{table}
\thispagestyle{empty}
\begin{tabular}{|ccccc|}
\hline $g$ & $n_g$ & $n_{g-1}+n_{g-2}$ & $\frac{n_{g-1}+n_{g-2}}{n_g}$ & $\frac{n_{g}}{n_{g-1}} $\\\hline 
0 & 1 & & & \\
1 & 1 & & & 1 \\
2 & 2 & 2 &
1
&
2
\\
3 & 4 & 3 &
0.75
&
2
\\
4 & 7 & 6 &
0.857143
&
1.75
\\
5 & 12 & 11 &
0.916667
&
1.71429
\\
6 & 23 & 19 &
0.826087
&
1.91667
\\
7 & 39 & 35 &
0.897436
&
1.69565
\\
8 & 67 & 62 &
0.925373
&
1.71795
\\
9 & 118 & 106 &
0.898305
&
1.76119
\\
10 & 204 & 185 &
0.906863
&
1.72881
\\
11 & 343 & 322 &
0.938776
&
1.68137
\\
12 & 592 & 547 &
0.923986
&
1.72595
\\
13 & 1001 & 935 &
0.934066
&
1.69088
\\
14 & 1693 & 1593 &
0.940933
&
1.69131
\\
15 & 2857 & 2694 &
0.942947
&
1.68754
\\
16 & 4806 & 4550 &
0.946733
&
1.68218
\\
17 & 8045 & 7663 &
0.952517
&
1.67395
\\
18 & 13467 & 12851 &
0.954259
&
1.67396
\\
19 & 22464 & 21512 &
0.957621
&
1.66808
\\
20 & 37396 & 35931 &
0.960825
&
1.66471
\\
21 & 62194 & 59860 &
0.962472
&
1.66312
\\
22 & 103246 & 99590 &
0.964589
&
1.66006
\\
23 & 170963 & 165440 &
0.967695
&
1.65588
\\
24 & 282828 & 274209 &
0.969526
&
1.65432
\\
25 & 467224 & 453791 &
0.971249
&
1.65197
\\
26 & 770832 & 750052 &
0.973042
&
1.64981
\\
27 & 1270267 & 1238056 &
0.974642
&
1.64792
\\
28 & 2091030 & 2041099 &
0.976121
&
1.64613
\\
29 & 3437839 & 3361297 &
0.977735
&
1.64409
\\
30 & 5646773 & 5528869 &
0.97912
&
1.64254
\\
31 & 9266788 & 9084612 &
0.980341
&
1.64108
\\
32 & 15195070 & 14913561 &
0.981474
&
1.63973
\\
33 & 24896206 & 24461858 &
0.982554
&
1.63844
\\
34 & 40761087 & 40091276 &
0.983567
&
1.63724
\\
35 & 66687201 & 65657293 &
0.984556
&
1.63605
\\
36 & 109032500 & 107448288 &
0.98547
&
1.63498
\\
37 & 178158289 & 175719701 &
0.986312
&
1.63399
\\
38 & 290939807 & 287190789 &
0.987114
&
1.63304
\\
39 & 474851445 & 469098096 &
0.987884
&
1.63213
\\
40 & 774614284 & 765791252 &
0.98861
&
1.63128
\\
41 & 1262992840 & 1249465729 &
0.98929
&
1.63048
\\
42 & 2058356522 & 2037607124 &
0.989919
&
1.62975
\\
43 & 3353191846 & 3321349362 &
0.990504
&
1.62906
\\
44 & 5460401576 & 5411548368 &
0.991053
&
1.62842
\\
45 & 8888486816 & 8813593422 &
0.991574
&
1.62781
\\
46 & 14463633648 & 14348888392 &
0.992067
&
1.62723
\\
47 & 23527845502 & 23352120464 &
0.992531
&
1.62669
\\
48 & 38260496374 & 37991479150 &
0.992969
&
1.62618
\\
49 & 62200036752 & 61788341876 &
0.993381
&
1.6257
\\
50 & 101090300128 & 100460533126 &
0.99377
&
1.62525
\\
\hline
\end{tabular}
\caption{Computational results on the number of numerical semigroups
  up to genus $\maxgenus$.}
\label{taula}
\end{table}

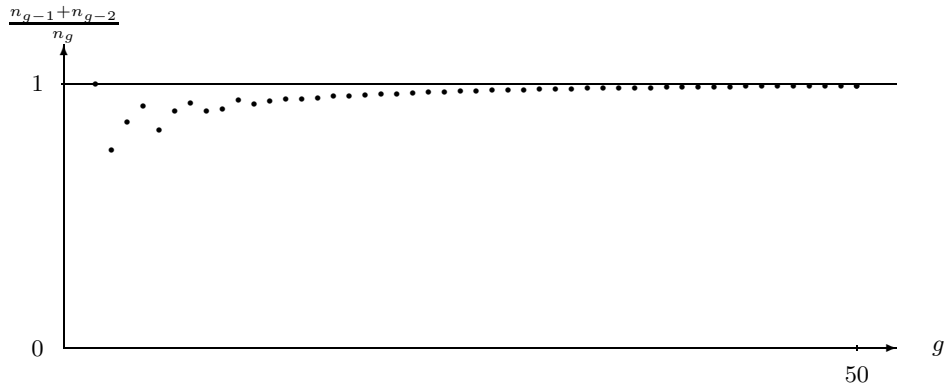
\begin{figure}
\begin{center}
\def\simbol{\circle*{2}}
\begin{picture}(350.000000,150.000000)
\put(20.000000,20.000000){\vector(1,0){315.000000}}
\put(20.000000,20.000000){\vector(0,1){115.000000}}
\put(342.500000,20.000000){\makebox(0,0){\phantom{mm}\small{$g$}}}
\put(20.000000,142.500000){\makebox(0,0){\small{$\frac{n_{g-1}+n_{g-2}}{n_g}$}}}
\put(10.000000,120.000000){\makebox(0,0){\small{1}}}
{\thinlines\put(20.000000,120.000000){\line(1,0){315.000000}}}
\put(10.000000,20.000000){\makebox(0,0){\small{0}}}
\put(19.000000,120.000000){\line(1,0){2.000000}}
\put(320.000000,10.000000){\makebox(0,0){\small{50}}}
\put(320.000000,19.000000){\line(0,1){2.000000}}
\put(32.000000,120.000000){\simbol}
\put(38.000000,95.000000){\simbol}
\put(44.000000,105.714298){\simbol}
\put(50.000000,111.666698){\simbol}
\put(56.000000,102.608700){\simbol}
\put(62.000000,109.743602){\simbol}
\put(68.000000,112.537302){\simbol}
\put(74.000000,109.830500){\simbol}
\put(80.000000,110.686297){\simbol}
\put(86.000000,113.877602){\simbol}
\put(92.000000,112.398602){\simbol}
\put(98.000000,113.406600){\simbol}
\put(104.000000,114.093299){\simbol}
\put(110.000000,114.294697){\simbol}
\put(116.000000,114.673300){\simbol}
\put(122.000000,115.251697){\simbol}
\put(128.000000,115.425898){\simbol}
\put(134.000000,115.762098){\simbol}
\put(140.000000,116.082503){\simbol}
\put(146.000000,116.247202){\simbol}
\put(152.000000,116.458900){\simbol}
\put(158.000000,116.769500){\simbol}
\put(164.000000,116.952599){\simbol}
\put(170.000000,117.124898){\simbol}
\put(176.000000,117.304201){\simbol}
\put(182.000000,117.464198){\simbol}
\put(188.000000,117.612101){\simbol}
\put(194.000000,117.773498){\simbol}
\put(200.000000,117.912002){\simbol}
\put(206.000000,118.034102){\simbol}
\put(212.000000,118.147398){\simbol}
\put(218.000000,118.255402){\simbol}
\put(224.000000,118.356700){\simbol}
\put(230.000000,118.455602){\simbol}
\put(236.000000,118.547000){\simbol}
\put(242.000000,118.631197){\simbol}
\put(248.000000,118.711401){\simbol}
\put(254.000000,118.788399){\simbol}
\put(260.000000,118.861003){\simbol}
\put(266.000000,118.929000){\simbol}
\put(272.000000,118.991901){\simbol}
\put(278.000000,119.050403){\simbol}
\put(284.000000,119.105299){\simbol}
\put(290.000000,119.157399){\simbol}
\put(296.000000,119.206698){\simbol}
\put(302.000000,119.253100){\simbol}
\put(308.000000,119.296898){\simbol}
\put(314.000000,119.338102){\simbol}
\put(320.000000,119.377000){\simbol}
\put(320.000000,119.377000){\simbol}
\end{picture}
\end{center}
\caption{Behavior of the quotient $\frac{n_{g-1}+n_{g-2}}{n_g}$.
The values in this graphic correspond to
the values in Table~\ref{taula}.}
\label{graficfib}
\end{figure}

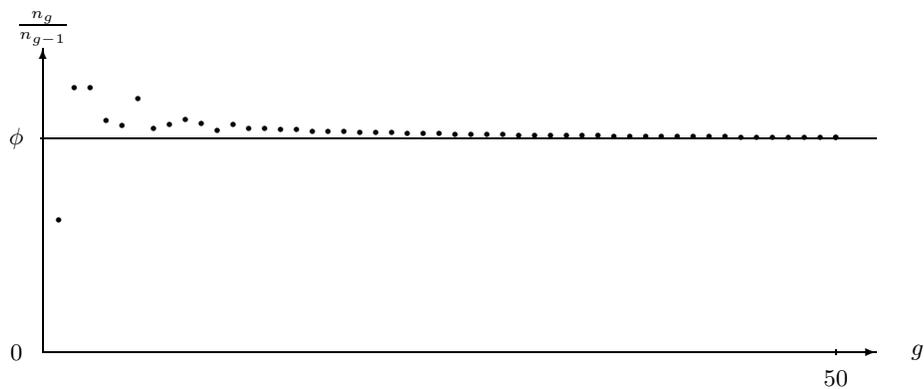
\begin{figure}
\begin{center}
\def\simbol{\circle*{2}}
\begin{picture}(350.000000,150.000000)
\put(20.000000,20.000000){\vector(1,0){315.000000}}
\put(20.000000,20.000000){\vector(0,1){115.000000}}
\put(342.500000,20.000000){\makebox(0,0){\phantom{mm}\small{$g$}}}
\put(20.000000,142.500000){\makebox(0,0){\small{$\frac{n_g}{n_{g-1}}$}}}
\put(10.000000,100.901699){\makebox(0,0){\small{$\phi$}}}
{\thinlines\put(20.000000,100.901699){\line(1,0){315.000000}}}
\put(10.000000,20.000000){\makebox(0,0){\small{0}}}
\put(19.000000,100.901699){\line(1,0){2.000000}}
\put(320.000000,10.000000){\makebox(0,0){\small{50}}}
\put(320.000000,19.000000){\line(0,1){2.000000}}
\put(26.000000,70.000000){\simbol}
\put(32.000000,120.000000){\simbol}
\put(38.000000,120.000000){\simbol}
\put(44.000000,107.500000){\simbol}
\put(50.000000,105.714501){\simbol}
\put(56.000000,115.833498){\simbol}
\put(62.000000,104.782499){\simbol}
\put(68.000000,105.897499){\simbol}
\put(74.000000,108.059503){\simbol}
\put(80.000000,106.440498){\simbol}
\put(86.000000,104.068501){\simbol}
\put(92.000000,106.297500){\simbol}
\put(98.000000,104.543997){\simbol}
\put(104.000000,104.565502){\simbol}
\put(110.000000,104.377003){\simbol}
\put(116.000000,104.109002){\simbol}
\put(122.000000,103.697498){\simbol}
\put(128.000000,103.697999){\simbol}
\put(134.000000,103.403999){\simbol}
\put(140.000000,103.235502){\simbol}
\put(146.000000,103.156002){\simbol}
\put(152.000000,103.003002){\simbol}
\put(158.000000,102.793999){\simbol}
\put(164.000000,102.716000){\simbol}
\put(170.000000,102.598501){\simbol}
\put(176.000000,102.490498){\simbol}
\put(182.000000,102.396001){\simbol}
\put(188.000000,102.306498){\simbol}
\put(194.000000,102.204503){\simbol}
\put(200.000000,102.126999){\simbol}
\put(206.000000,102.054001){\simbol}
\put(212.000000,101.986499){\simbol}
\put(218.000000,101.922001){\simbol}
\put(224.000000,101.862003){\simbol}
\put(230.000000,101.802499){\simbol}
\put(236.000000,101.748998){\simbol}
\put(242.000000,101.699502){\simbol}
\put(248.000000,101.651998){\simbol}
\put(254.000000,101.606501){\simbol}
\put(260.000000,101.563997){\simbol}
\put(266.000000,101.524003){\simbol}
\put(272.000000,101.487501){\simbol}
\put(278.000000,101.453001){\simbol}
\put(284.000000,101.421000){\simbol}
\put(290.000000,101.390500){\simbol}
\put(296.000000,101.361502){\simbol}
\put(302.000000,101.334502){\simbol}
\put(308.000000,101.309003){\simbol}
\put(314.000000,101.285000){\simbol}
\put(320.000000,101.262499){\simbol}
\put(320.000000,101.262499){\simbol}
\end{picture}
\end{center}
\caption{Behavior of the quotient $\frac{n_{g}}{n_{g-1}}$.
The values in this graphic correspond to
the values in Table~\ref{taula}.}
\label{graficor}
\end{figure}

In Table~\ref{taula} there are the results obtained for all numerical semigroups
with genus up to $\maxgenus$. For each genus we wrote the number of numerical semigroups of
the given genus, the Fibonacci-like-estimated value given by the sum of the number of semigroups
of the two previous genus, the value of the quotient
$\frac{n_{g-1}+n_{g-2}}{n_g}$,
and the value of the quotient $\frac{n_g}{n_{g-1}}$.
In Figure~\ref{graficfib} and Figure~\ref{graficor}
we depicted the behavior of these quotients. From
these graphics one can predict that
$\frac{n_{g-1}+n_{g-2}}{n_g}$
approaches $1$ as $g$ approaches infinity
whereas
$\frac{n_g}{n_{g-1}}$
approaches the golden ratio as $g$ approaches infinity.
We leave this as a conjecture.

\section{On the Wilf conjecture}
The Wilf conjecture (\cite{Wilf,DoMa})
states that the number $e$ of minimal generators of a numerical semigroup of
genus $g$ and conductor $c$ satisfies
$$e\geq\frac{c}{c-g}.$$
It is easy to check it when the numerical semigroup is symmetric, that
is, when $c=2g$.
In \cite{DoMa} the inequality is proved for many
other cases.
Here we proved by brute approach
that any numerical semigroup of genus at most $\maxgenus$
also satisfies the conjecture.

\section*{Acknowledgments}

The author would like to thank Jordi Funollet and Josep M. Mondelo
for their help on the computations and
Pedro A. Garc\'\i a-S\'anchez, Jos\'e
Carlos Rosales, Aureliano M. Robles, and Ruud Pellikaan for many valuable discussions.

This work was partly supported by the Spanish Ministry of
Science and Education through project CONSOLIDER INGENIO 2010
CSD2007-00004 ``ARES''.

\end{document}